\documentclass[10pt, english]{article}

\usepackage{latexsym, graphicx}
\usepackage{amsmath}
\usepackage{amsfonts}
\usepackage{amssymb}
\usepackage{hyperref}
\usepackage{babel}
\usepackage[babel=true]{csquotes}
\usepackage{geometry}

\usepackage{amssymb,amsmath,amsfonts,mathrsfs,amsxtra,amsthm}
\usepackage{graphicx,color,subfig}
\usepackage{accents,xspace}
\usepackage{dsfont}
\usepackage{amscd}
\newcommand{\mypapersize}{
\setlength{\textwidth }{16cm}
\setlength{\textheight}{23.5cm}
\setlength{\oddsidemargin}{-0.14cm}
\setlength{\topmargin}{-1.6cm}
}

\mypapersize

\newtheorem{thm}{Theorem}[section]
\newtheorem*{thm*}{Theorem}

\newtheorem{oss}[thm]{Remark}

\begin{document}
\title{Rebuttal\, to \vspace{0.2cm}\\
\enquote{Well-posedness and optimal decay rates for the viscoelastic Kirchhoff equation} \vspace{0.2cm}
}
\author{Fatiha Alabau-Boussouira\footnote{Laboratoire Jacques-Louis Lions Sorbonne Universit\'{e}, Universit\'{e} de Lorraine} 
}
\bibliographystyle{plain}

\date{March 31, 2022}
\maketitle {}

\begin{abstract}
A paper, entitled "Well-posedness and optimal decay rates for the viscoelastic Kirchhoff equation" was published in 2016, in the journal
{\it Boletim da Sociedade Paranaense de Matematica}. The main results of that paper are stated in Theorem 1.3. Section 2 is devoted to stability estimates for the solutions of the Galerkin approximation (2.1)-(2.2) of (1.2). Section 3 is devoted to the proof of the well-posedness for the Kirchhoff equation (1.2).
\  \\
This rebuttal is written in view to show the major flaws, and gaps in the proofs of this publication.
A major flaw stands in the proof of Lemma 2.2 related to well-known linear integral inequalities, so that the claimed stability estimates for the solutions of the Galerkin approximations (2.1)-(2.2) of (1.2) are not proved.
As a consequence, a major flaw also stands in Lemma 3.1 and in Theorem 1.3, the main result of the publication. 
The flaw in the proof of Lemma 2.2 appears in its conclusive part, which asserts, if one follows the claimed mathematical arguments, that any strictly positive, non-increasing function $H$ such that $\int_0^{\infty} H(t)dt < \infty$ should decay exponentially at infinity. This is false for an infinite number (uncountable) of such functions.
Counterexamples are given in the sequel.

In Section 3, the function $u$ involved in Lemma 3.1 and the sequel of the section are the solutions of (2.1)-(2.2), which is the Galerkin approximation of (1.2). A major concern is that nowhere the convergence of a subsequence of these solutions is considered, nor that of a proof that if the limit would exist, it would be a solution of the equation (1.2). Hence, the well-posedness of  (1.2) has not been established. A further concern is also that the publication does not contain any optimality results, even though such results are announced in the title, the abstract and in the introduction.

The publication also starts by an obvious contradictory assumption 1.1 which rises whenever one combines (1.3) and (1.4). One can assume for instance that $\gamma=0$ in (1.4) to avoid any contradiction, but in this case, the considered model is very close to the usual wave equation with viscoelastic dampings, and much stronger results are already available.
\end{abstract}

\vskip 0.2 cm

\section{Introduction}
We analyze the paper ”Well-posedness and optimal decay rates for the viscoelastic Kirchhoff equation”
published in the journal Boletim da Sociedade Paranaense de Matematica written by A\"issa Guesmia, Salim A. Messaoudi and Claudete M. Webler with the corresponding author A. Guesmia (see page 22).
The paper presents the partial integro-differential equation (1.1) introduced by Kirchhoff in the given reference [15].  This model depends on a nonnegative parameter $p_0$.
The authors mention that their paper handle both the degenerate model (when $p_0=0$) and the non-degenerate model (when $p_0>0$), and treat a more general model, namely the Kirchhoff-Carrier problem with memory, modeled by (1.2). The purpose of the paper is to \enquote{investigate the well-posedness as well as optimal decay rates estimates of the energy associated with} the so-called Kirchhoff-Carrier problem with memory, namely (1.2).

We present here only major flaws or gaps. Other small mistakes or typos are present in the paper. Since they can easily be corrected and do not affect in a major way the results, we do not detail them. 

Furthermore, we present the rebuttal in a pedagogic way for an unaware public of readers, so that they can easily understand the major flaws with concrete counterexamples, and the mathematical arguments explaining the flaws.

\vskip 0.2 cm 

A further important concern, is that there exists a second version compared to the present one analyzed in this first rebuttal, whose mathematical contents are, up to minor changes of indexation (in the equations or in the notation for the references, numbering of the main theorem) or presentation (equations on the left or right hand side), identical to the present published article. This second version is available in open access on ResearchGate.
 
 \vskip 0.2 cm 

 The title of this second version, after the first page of the ResearchGate presentation, is exactly the title of the paper published in the Boletim da Sociedade Paranaense de Matematica, and was written by \enquote{\small{M. M. CAVALCANTI$^{(1*)}$, V. N. DOMINGOS CAVALCANTI$^{(1)}$, A. GUESMIA$^{(2)}$ , S. A. MESSAOUDI$^{(3)}$, AND C. M. WEBLER${(1)}$}}

\vskip 0.1cm

This second version has thus five authors, the three present in the Boletim da Sociedade Paranaense de Matematica version, and ahead of the three authors two new authors :  M. M. Cavalcanti and V. N. Domingos Cavalcanti. M. M. Cavalcanti is presented as the corresponding author for this publication. 

\vskip 0.1cm

This second version further mentions that:
\vskip 0.2 cm
\enquote{Research of M. M. Cavalcanti partially supported by the CNPq Grant 300631/2003-0.Research of V. N. Domingos Cavalcanti partially supported by the CNPq Grant 304895/2003-2.}
\vskip 0.2 cm

\begin{oss}
On the website of the journal, one can notice that M. M. Cavalcanti and V. N. Domingos Cavalcanti, are respectively editor-in-chief and managing editor of the Boletim da Sociedade Paranaense de Matematica.
\end{oss}

\vskip 0.1cm

{\it The existence of two versions with similar mathematical contents (up to minor changes) is the reason why, this rebuttal has also a distinct second version, as a twin rebuttal, adapted to the second version with public access on ResearchGate, of the present published article in a peer-reviewed open access journal. Both rebuttals have been written in [4] and [5]}.

\section{A major flaw in the proof of Lemma 2.2 in Section 2. Counterexamples. Consequences on Section 3 and the main result}

The proof of Lemma 2.2 is finalized at pages 218-219 where the authors proved (2.47), that is
\vskip 0.2 cm
\enquote{
\begin{equation*}
\int_0^T \xi(t) \hat{E}(t) dt \leqslant C \hat{E}(0), \quad \forall \ T \geqslant t_0 \quad (2.47)
\end{equation*}
}
\vskip 0.2 cm
where the constant $C$ depends on $t_0$ and $g(t_0)$, and then at page 219 where the authors extend (2.47) to all $T>0$.
 Then the authors introduce 
 \vskip 0.2 cm
\enquote{
 $$
 \hat{\xi}(t) = \int_0^t \xi(s)ds \mbox{ and } F(t)= \hat{E}\left(  \hat{\xi}^{-1}(t)\right).
 $$
 }
 \vskip 0.2 cm
and then deduce that
\vskip 0.2 cm
\enquote{
$$
\int_0^T F(t) dt \leqslant CF(0), \quad \forall \ T \geqslant 0.
$$
}
\vskip 0.2 cm
Once this inequality is proved, the authors claim:
\vskip 0.3cm 
\enquote{Consequently, by applying Theorem 9.1 in} [3] (referenced as 16 in the paper analyzed here), enquote{we find that there exist positive constants $c$ and $\theta$ not depending on $\hat{E}(0)$ such that
$$
F(t) \leqslant cF(0) e^{-\theta \, t}, \quad \forall \ t \geqslant 0.
$$
}
\vskip 0.2cm

The authors then easily deduce the stability estimate (1.9) from the above inequality. The estimate (1.9) in the paper is presented in Theorem 1.3 at page 205 as:
\vskip 0.2 cm
\enquote{
\begin{equation*}
(1.9) \quad \quad \quad \hat{E}(t) \leqslant c \hat{E}(0)e^{-\theta\, \int_0^t \xi(s)ds}, \forall \quad t \geqslant 0,
\end{equation*}
where $\theta$ and $c$ are positive constants independent of the initial data.
}
\vskip 0.2 cm
Note that if the nonnegative, non-increasing function $F$, as involved above in the paper analyzed here, is such that there exists $T_0 \geqslant 0$ such that $F(T_0)=0$, then one has trivially that $F(t)\equiv 0$ on $[T_0,\infty)$. In particular, if $T_0=0$, $F$ is the null function.

Hence, the interesting cases are when $F$ is strictly positive on $[0,\infty)$.
Let us now discuss the above claim in the paper and show the major flaw.
\begin{oss}\label{Rk0}
Indeed, the required result in [3] (referenced as 16 in the paper analyzed here) is Theorem 8.1 and not Theorem 9.1 of [3] (see Remark~\ref{Rk3} below), that is
\begin{thm*}{\bf 8.1} (Reference [3])
\enquote{Let $E\ : \mathbb{R}_+ \longrightarrow \mathbb{R}_+ ( \mathbb{R}_+:= [0,\infty))$ be a non-increasing function and assume that there exists a constant $T>0$ such that
\begin{equation*}
\int_t^{\infty} E(s) ds \leqslant T E(t), \quad \forall \ t \in \mathbb{R}_+ .  
\end{equation*}
Then 
$$
E(t) \leqslant E(0) e^{1- \frac{t}{T}}, \quad \forall \ t \geqslant T.
$$
}
\end{thm*}
\end{oss}

A further observation in the reference [3] states that the above inequality also holds trivially for all $t \geqslant 0$.
\begin{oss}\label{Rk4}
Let us consider the function $F_2$ defined on $ \mathbb{R}_+$ by 
$$
F_2(t)=\dfrac{1}{(t+1)^2} \quad \forall \ t \geqslant 0.
$$
Then $F_2$ satisfies the inequality
$$
\int_0^T F_2(t) dt \leqslant F_2(0), \quad \forall \ T \geqslant 0, 
$$
whereas there does not exist any positive constants $c$ and $\theta$, such that $F_2$ satisfies the exponential decay announced in
$$
F_2(t) \leqslant cF_2(0) e^{-\theta \, t}, \quad \forall \ t \geqslant 0.
$$
\end{oss}
Now we can note more generally that 
\begin{oss}\label{Rk01}
\noindent The inequality 
$$
\int_0^T F(t) dt \leqslant CF(0), \quad \forall \ T \geqslant 0.
$$
proved in the paper analyzed here does not allow to apply the above Theorem 8.1 of [3] for the following mathematical arguments:

The inequality presented as an hypothesis in Theorem 8.1 of 16 ([3] here), states that any nonegative, nonincreasing function $F$ satisfying
\begin{equation*}
\int_t^{\infty} F(s) ds \leqslant C F(t), \quad \forall \ t \in \mathbb{R}_+ .  
\end{equation*}
decays exponentially at infinity. 
\vskip 0.2 cm
This inequality, implies in particular setting $t=0$ and using the nonnegativity of the function $F$ that 
\begin{equation*}
\int_0^{T} F(s) ds \leqslant C F(0), \quad \forall \ T \geqslant 0 .  
\end{equation*}
However the converse is false.
The major flaw in Lemma 2.2 of the paper stands in using this converse implication, which is false. 
\end{oss}
More generally, let us observe
\begin{oss}\label{Rk4bis}
If the conclusive part of the proof of Lemma 2.2 was true, it would imply that any strictly positive, non-increasing function $H$ such that $\int_0^{\infty} H(t)dt < \infty$ should decay exponentially at infinity. Let us explain why. Let $H$ be a strictly positive, non-increasing function such that $\int_0^{\infty} H(t)dt < \infty$. Then $H$ trivially satisfies
$$
\int_0^T H(t) dt \leqslant CH(0), \quad \forall \ T \geqslant 0,
$$
with 
$$
C=\dfrac{1}{H(0)} \,\int_0^{\infty} H(t) dt. 
$$
Then from the conclusive part of Lemma 2.2, $H$ should decay exponentially at infinity.
Let us explicit the following counterexamples (among others):

\vskip 0.2 cm

\noindent One can for instance define the function $H$ as follows

\vskip 0.2 cm

$$
H(t)= (t+1)^{-b} \quad \forall \ t \geqslant 0.
$$
where $b>1$ is arbitrary. It satisfies the estimate 
$$
\int_0^T H(t) dt \leqslant \dfrac{1}{b-1}H(0), \quad \forall \ T \geqslant 0,
$$
whereas it does not converge exponentially to zero at infinity. One can also consider functions $H$ of the form
$$
t \in [0,\infty) \longrightarrow H(t)=(t+e)^{-\alpha}\ln(t+e)^{-\beta},
$$
with either $\alpha>1$, or $\alpha =1$ and $\beta >1$ involved in the studies of Bertrand's integrals. This also leads to a series of counterexamples.

\end{oss}

\begin{oss}\label{Rk3}

\noindent Theorem 9.1 in the reference 16 ([3] here) of the paper concerns nonlinear integral inequalities, that is
\begin{thm*}{\bf 9.1} (Reference [3])
\enquote{Let $E\ : \mathbb{R}_+ \longrightarrow \mathbb{R}_+ ( \mathbb{R}_+:= [0,\infty))$ be a non-increasing function and assume that there exist two constants $\alpha >0$ and  $T>0$ such that
\begin{equation*}
\int_t^{\infty} E^{\alpha +1}(s) ds \leqslant T E(0)^{\alpha} E(t), \quad \forall \ t \in \mathbb{R}_+ .  
\end{equation*}
Then we have
$$
E(t) \leqslant E(0) \left(\dfrac{T + \alpha t}{T+ \alpha T}\right)^{-1/\alpha}, \quad \forall \ t \geqslant T.
$$
}
\end{thm*}
The nonlinear character of the integral inequality stands in the exponent $\alpha>0$. In the paper analyzed here, $\alpha=0$. Hence one needs to use a linear integral inequality result, i.e. Theorem 8.1 of [3].
\end{oss}

Hence, Lemma 2.2 at page 210 is not proved. Lemma 2.2 is used to get the inequality (3.7) in the proof of Lemma 3.1. As a consequence, Lemma 3.1 is also no longer valid. In a similar way, Theorem 1.3, the main result of the paper, is no longer valid.
\begin{oss}\label{Rk6}

\noindent Let us remark that even though Lemma 2.2 would have been proved, the estimate (1.9) at page 205 as part of the main result Theorem 1.3, would then have been proved only for the solutions of the approximate problem (2.2) for $m \in \mathbb{N}^{\ast}$ and not for the solution of the original problem (1.2). 
\end{oss}

\section{Missing proof of convergence of the Galerkin method and other missing arguments }
From the two above sections, we establish that the main result of the paper is not proved. Even though, let us assume that the previous steps are fine (even though they are not) and analyze the sequel of the proof.
\begin{oss}\label{Rk7}

\noindent In Lemma 3.1, $u$ has to be assumed to be a local solution of the approximate problem (2.2) on $[0,T[$ such that  the stated a priori bound holds. 
\end{oss}
We focus now on the text after the end of Remark 3.4 at page 221, that is on the text starting by 
\vskip 0.3 cm
\noindent \enquote{Recalling Lemma 2.2, one can assert that $u$ (the approximate solution constructed by the Galerkin method} until page 222 until the paragraph entitled \enquote{Uniqueness} :
\vskip 0.3 cm
The authors have chosen to suppress the subscript $m$ in this text whenever $u$, $u'$ are considered, for instance in (3.11),  (3.12), (3.13) and (3.14) for $u$ and $u'$ and also for $T$, which also has necessarily to be assumed as dependent of $m$ in (3.12), (3.13), and (3.14). This generates confusions in what is really proved and for which problem: the approximate problem (2.2) or the original problem (1.2). Indeed this is the approximate problem (2.2) which is considered in this part (see Remark~\ref{Rk9} below).

\vskip 0.1 cm

\begin{oss}\label{Rk8}

\noindent The authors consider initial data $(u_0,u_1) \in S_K$ where $K>0$ is fixed arbitrarily and want to prove (3.11) by contradiction for the solution of the approximate problem (2.2):

\vskip 0.3 cm

We first observe that since one deals with the approximate problem (2.2), one has to approximate the initial data $(u_0,u_1)$ as stated at page 207 in the two last lines of (2.2), where the approximated initial data $(u_{0m}, u_{1m})$ are assumed to converge in a stronger norm towards $(u_0,u_1)$, namely $\left(H^2(\Omega)\cap H^1_0(\Omega)\right)\times H^1_0(\Omega)$, as $m$ goes to $\infty$. It is further not specified how $\gamma_{jm}(0)$ and $\gamma_{jm}'(0)$ are chosen in (2.2). One may think that they correspond to the truncated series corresponding respectively to $u_0$ and $u_1$, obtained when the corresponding infinite sums are truncated to the sums until the integer $m$. 

\vskip 0.3 cm

One can easily check that the property that $(u_0,u_1) \in S_K$ as defined in (3.9) at page 221 does not necessarily imply that $(u_{0m}, u_{1m}) \in S_K$ for all $m \in \mathbb{N}^{\ast}$, even though one assumes that the initial data $(u_{0m}, u_{1m})$ are the truncated series until $m$ corresponding to $(u_0,u_1)$.
However, even with this choice, the property that $(u_{0m}, u_{1m}) \in S_K$ for all $m \in \mathbb{N}^{\ast}$ for all positive integers $m$ if $(u_0,u_1) \in S_K$, seems to hold true only  if $M$ is assumed to be non-decreasing, or if $m$ is sufficiently large. These arguments are missing at page 222 to apply the last inequality  in (3.14). Since no hypothesis of non-decreasingness on $M$ have been made, (3.14) holds true only  when $m$ is assumed to be sufficiently large, where the starting range at which it holds depends on the initial data $(u_0,u_1)$.
\end{oss}

The concluding part of the proof  of Lemma 3.1 is written as follows at page 222:
\vskip 0.3 cm
\noindent \enquote{As a consequence, we can repeat the continuation procedure indefinitely and we can conclude that, if $(u_0,u_1) \in S$, the solution $u$ can be continued globally on $\mathbb{R}_+$ and $(u(t),u'(t)) \in S$ for all $t \geqslant 0$.}
\vskip 0.3 cm
Let us observe that
\vskip 0.2 cm

\begin{oss}\label{Rk9} 
\begin{enumerate}
\medskip

\noindent \item It seems there is a confusion between the global solution $(u_m,u_m')$ of the approximate problem (2.2) considered at the beginning of the proof, and the continuous problem (1.2). The estimate (3.14) holds only for the approximate solutions and for $m \geqslant m_0$, where $m_0$ may depend on the initial data $(u_0,u_1)$. 

\item To be able to conclude, one needs to prove that the sequence of approximate solutions converge in a suitable norm, and that the limit solves problem (1.2) in a certain coherent sense. This would be the proof of the convergence of the Galerkin method, which is totally absent here.

\item Let us further remark that it does not seem the authors proved in this paper, or referred to a local existence of a solution $u$ for problem (1.2). This may holds true, but is not precised in the paper.
\end{enumerate}
\end{oss}
\section{Missing of results or discussion devoted to optimality}
We recall here all the parts of the article where the term \emph{optimal} is used, relatively to the quality of the claimed decay rates :
\begin{itemize}
\item In the title of the publication
\vskip 0.2 cm

\noindent \item In the abstract: (twice)
\vskip 0.2 cm

\noindent \item And further, twice in the introduction at page 204

\end{itemize}
\vskip 0.2 cm

Optimal decay rates are thus announced in the title, the abstract and the introduction, but nowhere else in the paper such matter is proved, nor discussed, especially in relation with the fact to deal with general decaying kernels $g$, involved in the memory term.

\section{Two contradictory assumptions in Section 1}

Two assumptions are introduced at page 205, before stating in the same page, the main result in Theorem 1.3. The first assumption is presented  as {\bf Assumption 1.1.} formulated as follows at the top of page 205:

\vskip 0.2 cm

\noindent \enquote{
{\bf Assumption 1.1.} 

\begin{eqnarray*}
\exists m_0>0 \, : M \in \mathcal{C}^1(\mathbb{R}_+) \mbox{ and }  M(\lambda) \geq m_0, & \forall \ \lambda \geq 0. & (1.3) \\
\exists \gamma\,, \delta>0 \, : M(\lambda) \leq \delta \lambda^{\gamma},& \forall \ \lambda \geq 0. & (1.4) \\
\exists \alpha\,, \beta>0 \, : \left|M'(\lambda)\right| \leq \beta \lambda^{\alpha}, & \forall \ \lambda \geq 0 . & (1.5)
\end{eqnarray*}
}

Let us observe
\begin{oss}\label{Rk1}
\begin{enumerate}
\noindent \item Since $m_0 >0$ in (1.3), together with $\gamma >0$ in (1.4), (1.3) and (1.4) are contradictory hypotheses and cannot hold together.  Indeed, setting $\lambda=0$ in both inequalities (1.3) and (1.4) leads to a contradiction. The assumption (1.3) corresponds to nondegenerate equations, whereas assumption (1.4) together with $M \geqslant 0$ corresponds to degenerate equations. A way to overcome these contradictions is to assume for instance $\gamma=0$ in (1.4), dealing then only with a nondegenerate equation. As mentioned in the abstract of this rebuttal, in this case, the considered model is very close to the usual wave equation with viscoelastic dampings, and much stronger results than the one claimed here, are already available.
\
\item The example (1.1) corresponds to $M(\lambda)= p_0  + \dfrac{Eh}{2L} \lambda, \quad \forall \  \lambda \geqslant 0$. It does not satisfies (1.3) unless $p_0>0$. In this latter case, it does not satisfy (1.4). In any case, it does not satisfy (1.5).
\end{enumerate}
\end{oss}
\begin{oss}\label{Rk2}
The inequality (1.3) of Assumption 1.1 is used at the following pages:

- page 208, last line

- page 211 in (2.15)

- page 217 in (2.40)

- page 220 in (3.5)

whereas the inequality (1.4) of Assumption 1.1 is used at the following pages:

- page 216 in (2.33)

- page 217 in (2.40).
\end{oss}

\section{References}

\bigskip

\noindent [1]  A. Guesmia, S. A. Messaoudi, C. M. Webler, Well-posedness and optimal decay rates for the viscoelastic Kirchhoff equation. Bol. Soc. Parana. Mat. 35 (3), 203--224, 2017, published online in 2016. DOI: https://doi.org/10.5269/bspm.v35i3.31395

\bigskip

\noindent [2] Cavalcanti-et-al-RG M. Cavalcanti, V. N. Domingos Cavalcanti, A. Guesmia, S. A. Messaoudi, C. M. Webler, Well-posedness and optimal decay rates for the viscoelastic Kirchhoff equation. Uploaded on ResearchGate by Marcelo M. Cavalcanti on 21 september 2018. \url{https://www.researchgate.net/publication/227598028_Existence_and_uniform_decay_for_nonlinear_viscoelastic_equation_with_strong_damping}

\bigskip

\noindent [3] V. Komornik, Exact Controllability and Stabilization. The Multiplier Method. RAM : Research in Applied Mathematics. Masson, Paris; John Wiley \& Sons, Ltd., Chichester, 1994.

\bigskip

\noindent [4] F. Alabau-Boussouira, Rebuttal\, of the publication \enquote{Well-posedness and optimal decay rates for the viscoelastic Kirchhoff equation}  in open access in the Boletim da Sociedade Paranaense de Matematica, Vol 35 (3) (2017).  \copyright{\ Fatiha Alabau-Boussouira in the original new Series: Mathema\&Ethics} (2021).

\bigskip

\noindent [5] F. Alabau-Boussouira, Rebuttal\, of the publication \enquote{Well-posedness and optimal decay rates for the viscoelastic Kirchhoff equation} published in open access on the ResearchGate Platform (2018). \copyright{\ Fatiha Alabau-Boussouira in the original new Series: Mathema\&Ethics} (2021).

\end{document}